\documentclass{amsart}

\theoremstyle{plain}
\newtheorem{thm}{Theorem}[section]

\theoremstyle{definition}
\newtheorem{defn}{Definition}[section]

\newtheorem{example}{Example}[section]

\numberwithin{equation}{section}

\begin{document}
\title{Homotonic Algebras}
\author{Michael Cwikel and Moshe Goldberg}

\address{Department of Mathematics,
Technion -- Israel Institute of Technology,
Haifa 32000, Israel}

\email{mcwikel@math.technion.ac.il \textit{and} goldberg@math.technion.ac.il}


\begin{abstract}
An algebra $\mathcal{A}$ of real or complex valued functions defined on a set $\mathbf{T}$ shall be called \textit{homotonic} if $\mathcal{A}$ is closed under forming of absolute values, and for all $f$ and $g$ in $\mathcal{A}$, the product $f\times g$ satisfies $|f\times g|\le|f|\times|g|$. Our main purpose in this paper is two-fold: To show that the above definition is equivalent to an earlier definition of homotonicity, and to provide a simple inequality which characterizes
sub-multiplicativity and strong stability for weighted sup norms on homotonic algebras.
\end{abstract}
\maketitle

\section{\label{sec:1}Definition and Examples}

Throughout this paper, let $\mathcal{A}$ denote a (finite or infinite dimensional) algebra over a field $\mathbb{F}$, either $\mathbb{R}$ or $\mathbb{C}$, of $\mathbb{F}$-valued functions defined on a given set $\mathbf{T}$. As usual, addition and scalar multiplication in $\mathcal{A}$ will be defined pointwise, i.e., for all $f$ and $g$ in $\mathcal{A}$, and all $\alpha$ in $\mathbb{F}$,
\begin{align*}
        &(f+g)(t)= f(t)+g(t),\\
        &(\alpha f)(t)= \alpha f(t).
    \end{align*}

Multiplication, often not pointwise, will be denoted by $\times$.

\begin{defn}
\label{def:1.1}
Let $\mathcal{A}$ be as above. We say that $\mathcal{A}$ is \textit{homotonic} if:

(i) $\mathcal{A}$ is closed under forming of absolute values, i.e., $f\in\mathcal{A}$ implies $|f|\in\mathcal{A}$.

(ii) For any two elements $f$ and $g$ in $\mathcal{A}$, we have $|f\times g|\le|f|\times|g|$.
\end{defn}

Here, for every $f\in\mathcal{A}$, the function $|f|$ is defined for each $t\in\mathbf{T}$ by $|f|(t)=|f(t)|$\,; and for real valued functions $f$ and $g$, the notation $f\le g$ will have the usual meaning, namely, $f(t)\le g(t)$ for all $t\in\mathbf{T}$.

We point out that Definition \ref{def:1.1} does not require $\mathcal{A}$ to be associative.

We also note that property (ii) implies that the product of non-negative functions in $\mathcal{A}$ is non-negative.

Examples of homotonic algebras are not hard to come by.
\begin{example}
\label{exa:1.1} Clearly, the algebra of all (bounded or not) $\mathbb{F}$-valued
functions defined on a given set $\mathbf{T}$, with pointwise multiplication
is homotonic.
\end{example}

\begin{example}
\label{exa:1.2} (Compare \cite{AG2}.) A more interesting example of a homotonic algebra is given by $\mathbb{F}^{n\times n}$, the algebra of all $n\times n$ matrices over $\mathbb{F}$ with the usual matrix operations. This algebra consists, of course, of all $\mathbb{F}$-valued functions on the set
\[
\mathbf{T}=\left\{ (j,k):\, j,k=1,...,n\right\} .
\]
\end{example}

\begin{example}
\label{exa:1.3} \cite{G}. To further illustrate homotonicity, fix positive constants $p$ and $\kappa$, and let $\mathcal{C}_{p,\kappa}(\mathbb{F})$ be the associative (and, in fact, commutative) algebra over $\mathbb{F}$ of all continuous, $p$-periodic, $\mathbb{F}$-valued functions on $\mathbb{R}$, where the product of $f$ and $g$ in $\mathcal{C}_{p,\kappa}(\mathbb{F})$ is defined by the convolution
\[
(f\ast g)(t)=\kappa\int_{0}^{p}f(t-x)g(x)dx,\quad t\in\mathbb{R}.
\]
Surely, if $f$ belongs to $\mathcal{C}_{p,\kappa}(\mathbb{F})$, so does $|f|$. Moreover, if $f$ and $g$ are members of $\mathcal{C}_{p,\kappa}(\mathbb{F})$, then
\begin{eqnarray*}
\left|f\ast g\right|(t) & = & \left|(f\ast g)(t)\right|
=\kappa\left|\int_{0}^{p}f(t-x)g(x)dx\right|\\
& \le & \kappa\int_{0}^{p}\left|f(t-x)g(x)\right|dx
=\kappa\int_{0}^{p}\left|f\right|(t-x)\left|g\right|(x)dx
=\left(\left|f\right|\ast\left|g\right|\right)(t),
\end{eqnarray*}
hence $\mathcal{C}_{p,\kappa}(\mathbb{F})$ is homotonic.

This example is a convenient prototype of many instances of algebras of functions defined on a locally compact abelian group where multiplication is a scalar multiple of convolution defined with respect to Haar measure on the group.
\end{example}

\begin{example}
\label{exa:1.4} Let $\mathcal{A}$ be a homotonic algebra, and let $\mathcal{A}^{+}$ be the algebra obtained by replacing the original product $f\times g$ in $\mathcal{A}$ by the Jordan product
\[
f\circ g\equiv\frac{1}{2}(f\times g+g\times f).
\]
Then it is not hard to see that $\mathcal{A}^{+}$ is also homotonic. Indeed, if $\mathcal{A}$ is closed under forming of absolute values, then so is $\mathcal{A}^{+}$. Further, if $f$ and $g$ are elements of $\mathcal{A}$ then, by the homotonicity of $\mathcal{A}$,
\[
|f\circ g|=\frac{1}{2}|f\times g+g\times f|\le\frac{1}{2}(|f\times g|+|g\times f|) \le\frac{1}{2}(|f|\times|g|+|g|\times|f|)=|f|\circ|g|.
\]
\end{example}

This example gives rise to straightforward constructions of \textit{non-associative} homotonic algebras. For instance, take $\mathcal{A}$ to be $\mathbb{F}^{n\times n}$ $(n\ge2)$, and consider $\mathbb{F}^{n\times n+}$, obtained by adopting the Jordan product $A\circ B\equiv\frac{1}{2}(AB+BA)$. For
\[
A=\left(\begin{array}{cc}
0 & 1\\
0 & 0\end{array}\right)\oplus O_{n-2}
,\quad B=\left(\begin{array}{cc}
0 & 0\\
1 & 0\end{array}\right)\oplus O_{n-2},
\]
$O_{n-2}$ denoting the ${(n-2)\times(n-2)}$ zero matrix, we have
\begin{eqnarray*}
(A\circ B)\circ B & = & \frac{1}{2}(AB+BA)\circ B
= \frac{1}{4}[(AB+BA)B+B(AB+BA)]\\
& = & \frac{1}{4}(AB^{2}+2BAB+B^{2}A)=\frac{1}{2}B
\end{eqnarray*}
and
\[
A\circ(B\circ B)=A\circ B^{2}=\frac{1}{2}(AB^{2}+B^{2}A)=0.
\]
Hence, $\mathbb{F}^{n\times n+}$ fails to be associative, although $\mathbb{F}^{n\times n}$ is.

\begin{example}
\label{exa:1.5}We note that if $\mathcal{B}$ is a sub-algebra of a homotonic algebra $\mathcal{A}$, then evidently, $\mathcal{B}$ is homotonic if and only if $\mathcal{B}$ is closed under forming of absolute values. For instance, consider the matrix algebra

\begin{equation}
\mathcal{A}_{2}(\mathbb{R})=\left\{ \left(\begin{array}{cc}
\alpha & \beta\\
-\beta & \alpha \end{array}
\right):\,\,\alpha,\beta\in\mathbb{R}\right\} \label{eq:1.1}
\end{equation}
with the usual matrix operations. Since this subalgebra of $\mathbb{R}^{2\times2}$ is not closed under forming of absolute values, it is not homotonic.
\end{example}

In the case where $\mathbb{F}=\mathbb{R}$, we can replace condition (ii) in Definition \ref{def:1.1} by a simpler condition:

\begin{thm}
\label{thm:1.1} Let $\mathcal{A}$ be an algebra over $\mathbb{R}$ of real valued functions defined on a given set $\mathbf{T}$. Then $\mathcal{A}$ is homotonic if and only if:

$\mathrm{(i)}$ $\mathcal{A}$ is closed under forming of absolute values, i.e.,
$f\in\mathcal{A}$ implies $|f|\in\mathcal{A}$.

$\mathrm{(ii)}_{\mathbb{R}}$ For each pair of non-negative functions $f$ and $g$ in $\mathcal{A}$, the product $f\times g$ is also non-negative.
\end{thm}

\noindent \textit{Proof.} If condition (ii) holds then, for any non-negative functions $f$ and $g$ in $\mathcal{A}$, we have
\[
f\times g=|f|\times|g|\ge|f\times g|\ge0,
\]
so (ii)$_{\mathbb{R}}$ is established. Hence, in order to complete the proof, it suffices to show that (i) and (ii)$_{\mathbb{R}}$ imply (ii).

Indeed, in view of (i), for each $u$ in $\mathcal{A}$ , the non-negative functions $u_{+}\equiv\frac{1}{2}(|u|+u)$ and $u_{-}\equiv\frac{1}{2}(|u|-u)$ are both in $\mathcal{A}$. Moreover, we have $u=u_{+}-u_{-}$ and $|u|=u_{+}+u_{-}$. Thus, for every $u$ and $v$ in $\mathcal{A}$,
\begin{equation}
|u|\times|v| = (u_{+}+u_{-})\times (v_{+}+v_{-}) = u_{+}\times v_{+}+u_{-}\times v_{-}
+u_{+}\times v_{-}+u_{-}\times v_{+}\,.\label{eq:1.2}
\end{equation}
and
\begin{equation}
u\times v = (u_{+}-u_{-})\times(v_{+}-v_{-}) = u_{+}\times v_{+}+u_{-}\times v_{-}
-u_{+}\times v_{-}-u_{-}\times v_{+}\,.\label{eq:1.3}
\end{equation}
By (ii)$_{\mathbb{R}}$, the products $u_{+}\times v_{+}$, $u_{-}\times v_{-}$,
$u_{+}\times v_{-}$ and $u_{-}\times v_{+}$ are all non-negative. So, comparing
(\ref{eq:1.2}) and (\ref{eq:1.3}), we get
\begin{eqnarray*}
\left|u\times v\right| & = & \left|u_{+}\times v_{+}+u_{-}\times v_{-}
-u_{+}\times v_{-}-u_{-}\times v_{+}\right|\\
& \le & u_{+}\times v_{+}+u_{-}\times v_{-}+u_{+}\times v_{-}
+u_{-}\times v_{+} = \left|u\right|\times\left|v\right|,
\end{eqnarray*}
and the proof follows. $\qed$
\begin{example}
\label{exa:1.6}To illustrate this theorem, consider the familiar real vector space
\[
\mathbb{R}^{2}=\{(\alpha,\beta):\alpha,\beta\in\mathbb{R}\}.
\]
For each $(\alpha,\beta)$ and $(\gamma,\delta)$ in $\mathbb{R}^{2}$, define multiplication by
\begin{equation}
(\alpha,\beta)\times(\gamma,\delta)
=(\alpha\gamma-\beta\delta,\alpha\delta+\beta\gamma)\,,\label{eq:1.4}
\end{equation}
which makes $\mathbb{R}^{2}$ into a $2$-dimensional algebra over the reals. Surely, $\mathbb{R}^{2}$ is closed under forming of absolute values, i.e., condition (i) holds. We observe, however, that if $\alpha$, $\beta$, $\gamma$ and $\delta$ are positive numbers with $\alpha\gamma<\beta\delta$, then the first component of the product $(\alpha,\beta)\times(\gamma,\delta)$
is negative; so condition (ii)$_{\mathbb{R}}$ fails, and by Theorem \ref{thm:1.1}, our algebra is not homotonic.

The mapping
\[
(\alpha,\beta)\mapsto\left(\begin{array}{cc}
\alpha & \beta\\
-\beta & \alpha\end{array}\right),\quad \alpha,\beta\in\mathbb{R},
\]
shows that the above algebra is an algebraically isomorphic image of the algebra $\mathcal{A}_{2}(\mathbb{R})$ defined in (\ref{eq:1.1}). In fact, the reader must have noticed by now that both these algebras are algebraically isomorphic to the complex numbers
\[
\mathbb{C}=\{\alpha+i\beta:\,\alpha,\beta\in\mathbb{R}\}
\]
viewed as a $2$-dimensional algebra over $\mathbb{R}$.
\end{example}

\section{\label{sec:2}An Earlier Equivalent Definition of Homotonic Algebras}

The notion of homotonicity was first introduced in \cite{AG2} in connection with functionals acting on a linear space $\mathbf{V}$ over $\mathbb{C}$ of bounded complex valued functions defined on a given set $\mathbf{T}$. In the same paper, the idea of homotonicity was extended to mappings from $\mathbf{V}$ into $\mathbf{V}$, and then to multiplication with which $\mathbf{V}$ was given the structure of an associative algebra.

Adapting the definitions in \cite{AG2}, the term \textit{homotonic algebra} was coined in \cite{G}. There, an associative algebra of bounded $\mathbb{F}$-valued functions defined on $\mathbf{T}$ is called \textit{homotonic} if:

(i) $\mathcal{A}$ is closed under forming of absolute values, i.e., $f\in\mathcal{A}$ implies $|f|\in\mathcal{A}$.

(ii)$^\prime$ For any four elements $f_{1}$, $f_{2}$, $g_{1}$ and $g_{2}$ in $\mathcal{A}$, such that $|f_{1}|\le g_{1}$ and $|f_{2}|\le g_{2}$, we have $|f_{1}\times f_{2}|\le g_{1}\times g_{2}$.

The name `homotonic' was chosen in this earlier definition because \textit{homo} indicates that multiplication preserves the relation $|f|\le g$ and \textit{tonic} reflects the fact that this relation is about order.

We shall now show that even in the general case where $\mathcal{A}$ is not necessarily associative and the functions in $\mathcal{A}$ are not necessarily bounded, the old and new definitions of homotonicity coincide. More precisely, we post:

\begin{thm}
\label{thm:2.1} Let $\mathcal{A}$ be an algebra over $\mathbb{F}$ of $\mathbb{F}$-valued functions defined on a set $\mathbf{T}$. Then $\mathcal{A}$ is homotonic if and only if conditions $\mathrm{(i)}$ and $\mathrm{(ii)^\prime}$ hold.
\end{thm}

\noindent \textit{Proof.} Putting $f_{1}=f$, $f_{2}=g$, $g_{1}=|f|$ and $g_{2}=|g|$, we immediately observe that (ii)$^\prime$ implies (ii). So assume that (i) and (ii) hold, and let us prove (ii)$^\prime$, thus forcing the desired result.

If $f$ and $g$ are non-negative functions in $\mathcal{A}$, then by (ii),
\[
f\times g=|f|\times|g|\ge|f\times g|\ge0;
\]
hence, as in the proof of Theorem \ref{thm:1.1}, (ii) implies (ii)$_{\mathbb{R}}$. Let $u$, $v$ and $w$ be real valued functions in $\mathcal{A}$ with $u\le v$ and $w\ge0$. Then, by (ii)$_{\mathbb{R}}$ ,
\[
v\times w-u\times w=(v-u)\times w\ge0;\]
so
\begin{equation}
u\le v\mbox{ and }w\ge0\Rightarrow u\times w\le v\times w\,.\label{eq:2.1}
\end{equation}
Analogously, we get
\begin{equation}
u\le v\,\mbox{and }w\ge0\Rightarrow w\times u\le w\times v\,.\label{eq:2.2}
\end{equation}

Suppose now that $f_{1}$, $f_{2}$, $g_{1}$ and $g_{2}$ are arbitrary functions in $\mathcal{A}$ which satisfy $|f_{1}|\le g_{1}$ and $|f_{2}|\le g_{2}$. Then appealing to (ii), (\ref{eq:2.1}) and (\ref{eq:2.2}) (in that order), we obtain
\[
|f_{1}\times f_{2}|\le|f_{1}|\times|f_{2}|\le g_{1}\times|f_{2}|\le g_{1}\times g_{2},
\]
and we are done.$\qed$

\section{\label{sec:3}Sub-multiplicative Weighted Sup Norms on Homotonic Algebras}

Our study of homotonic algebras is motivated mainly by the following theorem which provides a simple characterization of sub-multiplicativity for weighted sup norms.

Here, as usual, we call a norm on an algebra $\mathcal{A}$ \textit{sub-multiplicative} if
\[
\Vert f\times g\Vert\le\Vert f\Vert\Vert g\Vert\ \quad \mbox{for all }f,g\in\mathcal{A}.
\]

\begin{thm}
\label{thm:3.1} $\mathrm{(Compare}$ \cite[Theorem 4.2]{AG2}$\mathrm{.)}$ Let $\mathcal{A}$ be a homotonic algebra over $\mathbb{F}$ of $\mathbb{F}$-valued functions defined on a set $\mathbf{T}$. Let $w$ be a fixed positive function on $\mathbf{T}$ $\mathrm{(}$not necessarily in $\mathcal{A}  \mathrm{)}$, such that $w_{-1}$, defined by $w_{-1}(t)=\frac{1}{w(t)}$ for all $t\in\mathbf{T}$, is an element of $\mathcal{A}$. Assume that
\[
\sup_{t\in\mathbf{T}}w(t)|f(t)|<\infty\, \quad \mbox{for all }f\in\mathcal{A}.
\]
Then the weighted sup norm
\begin{equation}
\Vert f\Vert_{w,\infty}\equiv\sup_{t\in\mathbf{T}}w(t)|f(t)|, \quad f\in\mathcal{A},\label{eq:3.1}
\end{equation}
is sub-multiplicative on $\mathcal{A}$ if and only if
\begin{equation}
w_{-1}\times w_{-1}\le w_{-1}\,.\label{eq:3.2}
\end{equation}
\end{thm}

\noindent \textit{Proof.} Suppose that $\Vert\cdot\Vert_{w,\infty}$ is sub-multiplicative. Since $w_{-1}$ is a member of $\mathcal{A}\,$, it follows that
\begin{equation}
\Vert w_{-1}\times w_{-1}\Vert_{w,\infty}\le\Vert w_{-1}\Vert^{2}=1;\label{eq:3.3}
\end{equation}
hence,
\[
|w_{-1}\times w_{-1}|\le w_{-1}.
\]
Since $w_{-1}$ is a positive function, the homotonicity of $\mathcal{A}$ implies that $w_{-1}\times w_{-1}\ge0$; thus
\begin{equation}
w_{-1}\times w_{-1}=|w_{-1}\times w_{-1}|\le w_{-1} \label{eq:3.4}
\end{equation}
and (\ref{eq:3.2}) is in the bag.

Conversely, let (\ref{eq:3.2}) hold. Set
\[
\lambda\equiv\sup\{\Vert f\times g\Vert_{w,\infty}:~f,g\in\mathcal{A},~
\Vert f\Vert_{w,\infty}=\Vert g\Vert_{w,\infty}=1\},
\]
and observe that $\Vert\cdot\Vert_{w,\infty}$ is sub-multiplicative if and only if $\lambda\le1$. Select $f,g\in\mathcal{A}$ with $\Vert f\Vert_{w,\infty}=\Vert g\Vert_{w,\infty}=1$; hence \begin{equation}
|f|\le w_{-1}\,\mbox{and }|g|\le w_{-1}\,.\label{eq:3.5}
\end{equation}
Since $\mathcal{A}$ is homotonic, Theorem \ref{thm:2.1} guarantees that condition (ii)$^\prime$ holds. By (\ref{eq:3.5}), therefore,
\[
|f\times g|\le w_{-1}\times w_{-1},
\]
so aided by (\ref{eq:3.4}), we get
\[
|f\times g|\le w_{-1}.
\]
Consequently,
\[
\Vert f\times g\Vert_{w,\infty}\le\Vert w_{-1}\Vert_{w,\infty}=1;
\]
whence $\lambda\le1$, and the proof is complete.$\qed$

\begin{example}
\label{exa:3.1}(Compare \cite[Theorem 1]{AG1}.) To illustrate Theorem \ref{thm:3.1}, let us revisit $\mathbb{F}^{n\times n}$, the algebra of $n\times n$ matrices over $\mathbb{F}$ with the usual matrix operations. Let $W=(\omega_{jk})$ be a fixed $n\times n$ matrix of positive entries $\omega_{jk}$, and consider the weighted sup norm
\begin{equation}
\Vert A\Vert_{W,\infty}=\max_{j,k}\omega_{jk}|\alpha_{jk}|,
\quad A=(\alpha_{jk})\in\mathbb{F}^{n\times n}\,.\label{eq:3.6}
\end{equation}
Let $W_{-1}$ be the Hadamard inverse of $W$, that is, the matrix whose $(j,k)$ entry is $\frac{1}{\omega_{jk}}$. Then, by the theorem, $\Vert\cdot\Vert_{W,\infty}$ is sub-multiplicative if and only if
\begin{equation}
(W_{-1})^{2}\le W_{-1}\,,\label{eq:3.7}
\end{equation}
where $(W_{-1})^{2}$ is the usual squaring of $W_{-1}$, and where the inequality in (\ref{eq:3.7}) is construed entrywise. For instance (compare \cite[Corollary 1.1]{GS}), selecting $W=\mu E$, where $\mu$ is a positive constant and $E$ is the matrix all of whose entries are $1$, we easily find that the norm in (\ref{eq:3.6}) is multiplicative if and only if
\[
\mu\ge n.
\]
 In other words, the norm
\[
\Vert A\Vert_{\mu,\infty}\equiv\mu\max_{j,k}|\alpha_{jk}|,
\quad A=(\alpha_{jk})\in\mathbb{F}^{n\times n}\ ,
\]
is sub-multiplicative if and only if $\mu\ge n$.

Surely, the results in this example remain valid when the sup norm in (\ref{eq:3.6}) is applied to the \textit{non-associative} algebra $\mathbb{F}^{n\times n+}$ defined in Example \ref{exa:1.4}.
\end{example}

\begin{example}
\label{exa:3.2} \cite{G}. Falling back on the algebra $\mathcal{C}_{p,\kappa}(\mathbb{F})$ in Example \ref{exa:1.3}, we let $w$ be a continuous, $p$-periodic, positive function on $\mathbb{R}$. Then, evidently, $w_{-1}$ belongs to $\mathcal{C}_{p,\kappa}(\mathbb{F})$; so by Theorem \ref{thm:3.1}, the $w$-weighted sup norm
\[
\Vert f\Vert_{w,\infty}=\max_{0\le t\le p}w(t)|f(t)|,
\quad f\in\mathcal{C}_{p,\kappa}(\mathbb{F}),
\]
is sub-multiplicative if and only if $w_{-1}\ast w_{-1}\le w_{-1}$; that is, precisely when
\[
\kappa\int_{0}^{p}\frac{dx}{w(t-x)w(x)}\le\frac{1}{w(t)}, \quad 0\le t\le p.
\]
In particular, we see that the usual sup norm
\[
\Vert f\Vert_{\infty}=\max_{0\le t\le p}|f(t)|, \quad f\in\mathcal{C}_{p,\kappa}(\mathbb{F}),
\]
is sub-multiplicative if and only if $\kappa p\le1$.
\end{example}

Our next example involves an algebra of \textit{unbounded} functions where the weight function $w$ is \textit{not} a member of $\mathcal{A}$.

\begin{example}
\label{exa:3.3} Set $\mathbf{T}=(0,\infty)$, and let $\mathcal{A}$ be the real vector space of all functions on $\mathbf{T}$ of the form $f(t)=\alpha t$ where $\alpha$ is a real constant. For each
$f$ and $g$ in $\mathcal{A}$, define the product $f\times g$ by
\[
(f\times g)(t)=\frac{f(t)g(t)}{t}, \quad t\in\mathbf{T},
\]
thus making $\mathcal{A}$ into a homotonic algebra which is a faithful image of $\mathbb{R}$. Let $w:\mathbf{T}\to\mathbb{R}$ be the positive unbounded function $w(t)=\nu t^{-1}$ where $\nu$ is a positive constant. Note that $w$ is not an element of $\mathcal{A}$ but $w_{-1}$ is. With this choice of $w$, and for each $f(t)=\alpha t$ in $\mathcal{A}$, we have
\[
\sup_{t\in\mathbf{T}}w(t)|f(t)|=\nu|\alpha|<\infty\ .
\]
Hence, by Theorem \ref{thm:3.1}, the weighted sup norm
\begin{equation}
\Vert f\Vert_{w,\infty}=\sup_{t\in\mathbf{T}}w(t)|f(t)|\label{eq:3.8}
\end{equation}
is sub-multiplicative on $\mathcal{A}$ if and only if $w_{-1}\times w_{-1}\le w_{-1}$; that is, if and only if $\nu\ge1$.
\end{example}

\section{\label{sec:4} Strongly Stable Weighted Sup Norms on Homotonic Algebras}

As usual, whether the algebra $\mathcal{A}$ is associative or not, we define powers of each element $f\in\mathcal{A}$ inductively by
\begin{align*}
&f^1 = f,\\
&f^k = f^{k-1}\times f, \quad k=2,3,4,...
\end{align*}

Having powers at our disposal, we follow standard nomenclature and say that a norm $\Vert\cdot\Vert$ on $\mathcal{A}$ is \textit{strongly stable} if
\[
\Vert f^{k}\Vert\le\Vert f\Vert^{k} \quad \mbox{for all }f\in\mathcal{A}\ \mbox{and }k=1,2,3,...
\]

With these definition, we can now easily characterize strong stability for weighted sup norms on  homotonic algebras.

\begin{thm}
\label{thm:4.1} $\mathrm{(Compare}$ \cite[Theorem 4.2]{AG2}$\mathrm{.)}$ Let $\mathcal{A}$ be a homotonic algebra over $\mathbb{F}$ of $\mathbb{F}$-valued functions defined on a set $\mathbf{T}$. Let $w$ be a fixed positive function on $\mathbf{T}$ $\mathrm{(}$not necessarily in $\mathcal{A} \mathrm{)}$, such that $w_{-1}$ belongs to $\mathcal{A}$ . Assume that
\[
\sup_{t\in\mathbf{T}}w(t)|f(t)|<\infty\ \quad \mbox{for all }f\in\mathcal{A}.
\]
Then the weighted sup norm $\Vert\cdot\Vert_{w,\infty}$ in $\mathrm{(\ref{eq:3.1})}$ is strongly stable if and only if
\[
w_{-1}\times w_{-1}\le w_{-1}.
\]
\end{thm}

\noindent \textit{Proof.} If $w_{-1}\times w_{-1}\le w_{-1}$ holds, then by Theorem \ref{thm:3.1}, $\Vert\cdot\Vert_{w,\infty}$ is sub-multiplicative, hence strongly stable since for all $f$ in $\mathcal{A}$,
\[
\Vert f^{k}\Vert = \Vert f^{k-1}\times f\Vert \le\Vert f^{k-1}\Vert \Vert f\Vert, \quad k=2,3,4...
\]
Conversely, if $\Vert\cdot\Vert_{w,\infty}$ is strongly stable, then
\[
\Vert f\times f\Vert_{w,\infty}\le\Vert f\Vert_{w,\infty}^{2} \quad \mbox{for all } f\in\mathcal{A}.
\]
So setting $f=w_{-1}$, we get (\ref{eq:3.3}) and $w_{-1}\times w_{-1}\le w_{-1}$ follows.$\qed$

\medskip

Theorems 3.1 and 4.1 show, of course, that in the homotonic case, sub-multiplica\-tivity and strong stability are equivalent for weighted sup norms. It thus follows that the examples presented in Section \ref{sec:3} are also relevant here, in the sense that in each of those examples, the condition given for sub-multiplicativity is also necessary and sufficient for strong stability.

We conclude by remarking that in general, a strongly stable norm on a (homotonic) algebra may fail to be sub-multiplicative. A familiar example is the numerical radius,
\[
r(A)=\max\{|(Ax,x)|:~x\in\mathbb{C}^n,~(x,x)=1\},
\]
defined on $\mathbb{C}^{n\times n}$ $(n\ge 2)$, with respect to a given inner product $(\cdot ,\cdot)$ on $\mathbb{C}^n$. It is well known (e.g., \cite[Chapter 17]{H}) that $r$ is a norm on $\mathbb{C}^{n\times n}$ which is not sub-multiplicative; on the other hand, the celebrated Berger Inequality, \cite{B,P}, tells us that $r$ is strongly stable.

\end{document}